\documentclass{nseJournal}

\usepackage{amsfonts}

\renewcommand{\vec}[1]{\boldsymbol{#1}}
\newcommand{\bx}{\vec{x}}
\newcommand{\bOmega}{\vec{\Omega}}
\newcommand{\bnormal}{\vec{\hat{n}}}

\newcommand{\bgrad}{\vec{\nabla}}
\newcommand{\bdiv}{\bgrad \cdot}

\newcommand{\half}{\dfrac{1}{2}}
\newcommand{\third}{\dfrac{1}{3}}

\newcommand{\bpsi}{\vec{\psi}}
\newcommand{\bphi}{\vec{\phi}}
\newcommand{\bvarphi}{\vec{\varphi}}
\newcommand{\bJ}{\vec{J}}
\newcommand{\bP}{\vec{P}}
\newcommand{\bT}{\vec{T}}

\newcommand{\jump}[1]{\left[\!\left[ #1 \right]\!\right]}
\newcommand{\average}[1]{\left\{\!\!\left\{ #1 \right\}\!\!\right\}}

\begin{document}

\title{A Second Moment Method for $k$-Eigenvalue Acceleration with Continuous Diffusion and Discontinuous Transport Discretizations} %title of paper

% Use the \addAuthor macro to add authors in the order they should appear. The second argument corresponds to
% the affiliation declared below.
% The corresponding author should be wrapped in \correspondingAuthor
\addAuthor{\correspondingAuthor{Zachary K. Hardy}}{a}

% The corresponding author's email can be specified using \correspondingEmail
\correspondingEmail{zhardy@lanl.gov}
\addAuthor{Jim E. Morel}{b}
\addAuthor{Jan I.C. Vermaak}{c}
% Affiliations can be added in the order they should appear. For breaks in addresses, use either \\ or \tabularnewline
% \addAffiliation{a}{Sandia National Laboratories, Radiation Effects Theory Department\\ P.O. Box 5800, MS1179, Albuquerque, New Mexico 87185}
\addAffiliation{a}{Los Alamos National Laboratory, X-Theoretical Design Division\\ P.O. Box 1663, Los Alamos, New Mexico, 87545}
\addAffiliation{b}{Texas A\&M University, Department of Nuclear Engineering\\ AI Engineering Building\\ 423 Spence Street, MS 3133, College Station, Texas 77843-3133}
\addAffiliation{c}{Idaho National Laboratory, Reactor Physics Methods and Analysis Department\\ Idaho Falls, Idaho, 83402}

% Add keywords to appear in Abstract in the order they should appear
\addKeyword{radiation transport}
\addKeyword{$k$-eigenvalue}
\addKeyword{acceleration}
\addKeyword{second moment method}

\titlePage

\begin{abstract}

The second moment method is a linear acceleration technique that couples the transport equation to a diffusion equation with transport-dependent additive closures. 
The resulting low-order diffusion equation can be discretized independent of the transport discretization, unlike diffusion synthetic acceleration, and is symmetric positive definite, unlike quasidiffusion.
While this method has been shown to be comparable to quasidiffusion in iterative performance for fixed source and time-dependent problems, it is largely unexplored as an eigenvalue problem acceleration scheme due to the thought that the resulting inhomogeneous source makes the problem ill posed.
Recently, a preliminary feasibility study was performed on the second moment method for eigenvalue problems.
The results suggested comparable performance to quasidiffusion and more robust performance than diffusion synthetic acceleration.
This work extends the initial study to more realistic reactor problems using state-of-the-art discretization techniques.
Results in this paper show that the second moment method is more computationally efficient than its alternatives on complex reactor problems with unstructured meshes.

\end{abstract}

\section{Introduction}

The second moment method (SMM) was first introduced by Lewis and Miller in 1976 \cite{lewis1976comparison} for accelerating iterative solutions to the transport equation.
The method iteratively couples the transport equation to a diffusion equation with discrete, transport-dependent additive closures.
With the additive closures, the diffusion system is merely a reformulation of the transport equation, and as a result, captures relevant transport physics \cite{adams2002fast}.
This method belongs to the class of methods known as high-order/low-order (HOLO) methods \cite{chacon2017multiscale}.
In general, HOLO methods use the high-order (HO) system only to update closures in the equivalent low-order (LO) system, which converges the slowly convergent physics.
As a result of this, these methods generally have a reduced memory footprint, are rapidly convergent, and are more easily coupled to other physics in multiphysics simulations.

The SMM is closely related to a number of different acceleration schemes.
Via a formal linearization procedure, the SMM can be derived from quasidiffusion (QD) \cite{goldin1964quasi}, also known as the variable Eddington factor (VEF) method \cite{adams2002fast,olivier2023smm}.
The formulation of the two methods varies only in the choice of closures for the LO system.
While the SMM uses an additive closure to obtain a source corrected diffusion equation, QD uses multiplicative closures, resulting in a drift-diffusion equation \cite{olivier2023smm}.
Additionally, Tutt et al. \cite{tutt2023second} and Woodsford et al. \cite{woodsford2023variant} recently showed that the SMM is analytically equivalent to Alcouffe's source correction formulation of linear diffusion synthetic acceleration (DSA).

Each of these methods have their own advantages and disadvantages.
For example, because DSA methods solve a diffusion equation, the resulting systems are generally symmetric positive definite (SPD), allowing for the use of highly efficient preconditioned linear solvers, such as algebraic multigrid (AMG) preconditioned \cite{ruge1987algebraic} conjugate gradient (CG).
As is well studied, however, DSA methods require that the discrete diffusion equation be algebraically consistent with the discrete transport equation for unconditional stability \cite{reed1971effectiveness}.
While consistent diffusion discretizations have been developed for a number of different transport discretizations \cite{alcouffe1977diffusion,warsa2002fully,adams1992diffusion,wang2010diffusion}, they are considerably more complex and generally have more unknowns than standard diffusion discretizations.
On the other hand, QD can be discretized independent of the transport discretization with more familiar techniques.
When discretized independently, it should be noted that the QD solution differs from the transport solution on the order of the spatial discretization error \cite{olivier2023smm,olivier2023vef}.
Because QD results in a drift-diffusion equation, however, the resulting matrix is generally not SPD.
As a result of this, the aforementioned preconditioned linear solvers cannot be used \cite{olivier2023vef}.
The SMM inherits the favorable qualities of both DSA and QD without their respective drawbacks.
It not only can be discretized independently, but also yields a SPD system.
With the Fourier analysis done by Cefus and Larsen \cite{cefus1989stability} showing that QD and the SMM have similar iterative convergence properties, the SMM offers an attractive alternative to both DSA and QD.

Historically, nonlinear acceleration techniques for $k$-eigenvalue problems have been the standard.
This is largely due to the belief that the inclusion of source corrections required by linear methods results in the eigenvalue problem being ill posed.
In fact, Alcouffe introduced two nonlinear $k$-eigenvalue DSA acceleration schemes in the aforementioned paper that introduced the source correction form of DSA for fixed source problems \cite{alcouffe1977diffusion}.
Today, nonlinear methods like coarse mesh finite difference (CMFD) \cite{smith2002full} and nonlinear diffusion acceleration (NDA) \cite{willert2014comparison} are widely used in the reactor physics community for $k$-eigenvalue problem acceleration.
These methods inherit some of the limitations of both QD and DSA. For example, they both require solutions to drift-diffusion equations, while neither have unconditional stability independent of the spatial discretization.

A small number of linear acceleration schemes for $k$-eigenvalue problems have been explored in the past.
In the 1980s Gelbard et al. developed a method that utilized the Fredholm alternative theorem (FAT) \cite{ramm2001simple}, requiring forward and adjoint diffusion $k$-eigenvalue problem solutions \cite{gelbard1982acceleration}.
More recently, in 2019 a similar algorithm named semilinear diffusion acceleration (SDA), which also relied on the FAT, was developed by Dodson et al. \cite{dodson2019sda}.
In that same year, the Barbu-Adams method, which features a diffusion $k$-eigenvalue problem with an inhomogeneous source term for an additive correction, was introduced \cite{barbu2019linear}.
In a more recent publication, Barbu and Adams showed that the Barbu-Adams method, when viewed as a whole, is a well-posed eigenvalue problem that does not require the use of the FAT \cite{barbu2023convergence}.

In Tutt et al. \cite{tutt2023second} and Woodsford et al. \cite{woodsford2023variant}, it was shown that a SMM $k$-eigenvalue acceleration scheme is analytically equivalent to the Barbu-Adams method \cite{barbu2019linear}, but exhibits more favorable properties.
Namely, the SMM scheme maintains the SPD system obtained by the Barbu-Adams method, but allows for independent discretization of the LO system.
In Tutt et al. \cite{tutt2023second} and Woodsford et al. \cite{woodsford2023variant}, this method was demonstrated using a lumped linear discontinuous (LLD) discretization for transport and cell-centered diffusion discretizations.
The iterative performance of the method was almost identical to a QD $k$-eigenvalue acceleration scheme and comparable to the Barbu-Adams method.
In cases with optically thick cells, however, the Barbu-Adams method diverged while QD and the SMM did not.

The aim of this paper is to extend the initial work of Tutt et al. \cite{tutt2023second} and Woodsford et al. \cite{woodsford2023variant} to multiple spatial dimensions, unstructured meshes, more advanced spatial discretizations, and multiple energy groups.
In particular, the piece-wise linear discontinuous (PWLD) finite element discretization of Bailey and Adams \cite{bailey2008piecewise} is used to discretize transport along with a piece-wise linear continuous (PWLC) finite element discretization for diffusion.
The method is compared primarily to the Barbu-Adams method and a nonlinear formulation of the eigenvalue problem that uses a Jacobian-free Newton Krylov (JFNK) solver. 
For the prior, a PWLD transport discretization is used alongside the consistent modified interior penalty (MIP) PWLD diffusion discretization of Wang and Ragusa \cite{wang2010diffusion}.
In addition to showcasing the ability of the method to accelerate problems with high dominance ratios, the method is demonstrated in parallel on a realistic reactor geometry meshed with unstructured elements.

The remainder of this paper is organized as follows. 
In Section \ref{sec:method}, the SMM for $k$-eigenvalue problems is defined in a general context. 
This sections details the derivation of the multi-group moment equations with isotropic scattering, its boundary conditions, and the closures. 
In Section \ref{sec:discretization}, the discretization of the moment equations is discussed.
Following this, Section \ref{sec:results} presents results for two example problems.
First, a simple two-region, two-group slab problem with high dominance and scattering ratios is considered.
Here, the iterative performance of the SMM is compared to both standard power iterations (PI) and other acceleration schemes.
Additionally the order of convergence of the method SMM is estimated.
Next, the method is demonstrated on the C5G7 benchmark problem \cite{lewis2001c5g7}, where the SMM is compared against other acceleration schemes using an efficiency metric for both a varying number of discrete directions and a varying number of parallel processes.
Lastly, Section \ref{sec:conclusion} presents concluding remarks and potential extensions of this work.

\section{The Second Moment Method for Eigenvalue Problems}
\label{sec:method}

Consider the $k$-eigenvalue transport problem with isotropic scattering given by
\begin{equation}
    \begin{aligned}
        \bOmega \cdot &\bgrad \psi(\bx, E, \bOmega) + \sigma_t(\bx, E) \psi(\bx, E, \bOmega) = \\ &\dfrac{1}{4\pi}\left[ \int_{0}^{\infty} \sigma_s(\bx, E' \rightarrow E) \phi(
        \bx, E') dE' + \dfrac{\chi(\bx, E)}{k} \int_{0}^{\infty} \nu\sigma_f(\bx, E') \phi(\bx, E') dE' \right],
        \label{eq:lbs}
    \end{aligned}
\end{equation}
with vacuum and reflective boundary conditions
\begin{subequations}
    \begin{alignat}{2}
        \psi(\bx, E, \bOmega) &= 0, \qquad &&\bx \in \partial\Gamma_\text{vac}, \quad \bOmega \cdot \bnormal < 0 \\
        \psi(\bx, E, \bOmega) &= \psi(\bx, E, \bOmega'), \qquad &&\bx \in \partial\Gamma_\text{refl}, \quad \bOmega' = \bOmega - 2(\bOmega \cdot \bnormal) \bnormal,
    \end{alignat}
\end{subequations}
where $\bx \in \mathbb{R}^d$ is the spatial coordinate in dimension $d$ defined within domain $\Gamma \subset \mathbb{R}^d$ bounded by $\partial\Gamma \subset \mathbb{R}^{d-1}$, $E \in (0, \infty)$ is the particle energy, and $\bOmega \in \mathbb{S}^2$ is the particle direction on the unit sphere. 
The angular flux is given by $\psi$ and the scalar flux $\phi$ is defined by $\phi = \int_{4\pi} \psi d\bOmega$. 
The total cross section is denoted $\sigma_t$, the differential scattering cross section $\sigma_s$, and the fission cross section $\sigma_f$. $\chi$ is the fission spectrum and $\nu$ the average neutron yield per fission. 
The notation $E' \rightarrow E$ within the differential scattering cross section denotes scattering from energy $E'$ to energy $E$.

Applying the standard multi-group discretization in energy and discrete ordinates (DO) discretization in angle to Eq. \eqref{eq:lbs}, one obtains 
\begin{subequations}
    \begin{equation}
        \begin{aligned}
            \bOmega_n \cdot \bgrad &\psi_{ng}(\bx) + \sigma_{t,g}(\bx) \psi_{ng}(\bx) = \\&\dfrac{1}{4\pi} \left[ \sum_{g'=1}^{G} \sigma_{s,g'g}(\bx) \phi_{g'}(\bx) + \dfrac{\chi_g(\bx)}{k} \sum_{g'=1}^{G} \nu\sigma_{f,g'}(\bx) \phi_{g'}(\bx)  \right],
        \end{aligned}
        \label{eq:mgdo}
    \end{equation}\vspace{-8ex}
    \begin{alignat}{3}
    &\psi_{ng}(\bx) = 0, \qquad &&\bx \in \partial\Gamma_\text{vac}, \quad &&\bOmega_n \cdot \bnormal < 0, 
    \\
    &\psi_{ng}(\bx) = \psi_{n'g}(\bx), \qquad &&\bx \in \partial\Gamma_\text{refl}, \quad &&\bOmega_{n'} = \bOmega_n - 2(\bOmega_n \cdot \bnormal) \bnormal,
    \end{alignat}
\end{subequations}
where each equation is defined for energy groups $g = 1, \ldots, G$, and discrete directions $n = 1,\ldots, N$ \cite{lewis1984computational}. 
Here, $G$ and $N$ are the total number of energy groups and discrete directions, respectively.
The discrete directions are chosen such that they form a suitable angular quadrature from which flux moments can be computed.
The spatial discretization is left unspecified.

\subsection{Derivation of the Moment Equations}

The derivation of the moment equations begins by taking the zeroth and first angular moments of Eq. \eqref{eq:mgdo}
\begin{subequations}
    \begin{gather}
        \bdiv \bJ_g + \sigma_{r,g} \varphi_g = \sum_{g' \ne g} \sigma_{s,g'g} \varphi_{g'} + \dfrac{\chi_g}{k} \sum_{g'=1}^{G} \nu\sigma_{f,g'} \varphi_{g'},
        \label{eq:mg_zeroth}\\
        \bdiv \bP_g + \sigma_{t,g} \bJ_g = \vec{0}, 
        \label{eq:mg_first}
    \end{gather}
    \label{eq:mg_moment_system}
\end{subequations}
where all arguments are dropped for clarity and $\varphi$ is used in lieu of $\phi$ to distinguish between the transport and moment system scalar fluxes.
Here, $\sigma_{r,g} = \sigma_{t,g} - \sigma_{s,gg}$ is the group removal cross section.
The scalar flux (zeroth moment) $\varphi \in \mathbb{R}^1$, current (first moment) $\bJ \in \mathbb{R}^d$, and pressure (second moment) $\bP \in \mathbb{R}^{d \times d}$ are given by
\begin{equation*}
    \varphi = \sum_{n=1}^{N} w_n \psi_n, \qquad 
    \bJ = \sum_{n=1}^{N} w_n \bOmega_n \psi_n , \qquad
    \bP = \sum_{n=1}^{N} w_n (\bOmega_n \otimes \bOmega_n) \psi_n,
\end{equation*}
respectively, where the group index $g$ is dropped for clarity.
It is important to emphasize that the moments in Eq. \eqref{eq:mg_moment_system} are unknowns of the moment system, and not known functions of the angular flux.
Solving for the current in Eq. \eqref{eq:mg_first} and plugging the result into \eqref{eq:mg_zeroth} then gives
\begin{equation}
    -\bdiv \dfrac{1}{\sigma_{t,g}} \bdiv \bP_g + \sigma_{r,g} \varphi_g = \sum_{g' \ne g} \sigma_{s,g'g} \varphi_{g'} + \dfrac{\chi_g}{k} \sum_{g'=1}^{G} \nu\sigma_{f,g'} \varphi_{g'}.
    \label{eq:mg_moment}
\end{equation}

The boundary conditions for the moment equations are obtained by preserving the particle inflow at the boundary. 
Following from  Miften and Larsen \cite{miften1993quasi} and Olivier and Haut \cite{olivier2023smm}, and Olivier et al. \cite{olivier2023vef}, a vacuum boundary condition can be written as
    \begin{equation}
    \bJ \cdot \bnormal = \sum_{n=1}^{N} w_n \lvert \bOmega_n \cdot \bnormal \rvert \psi_n = B, \qquad \bx \in \partial\Gamma_\text{vac},
\end{equation}
where $B$ is an unknown boundary factor.
Reflective boundary conditions imply a zero net current, and are therefore trivially defined with
\begin{equation}
    \bJ \cdot \bnormal = 0, \qquad \bx \in \partial\Gamma_\text{refl}.
\end{equation}
With this, the moment equations are fully specified.

Due to the streaming operator, there are always more angular moment unknowns than equations in the moment equations. 
After eliminating the current in terms of the pressure in Eq. \eqref{eq:mg_moment}, there are $d^2 + 1$ unknowns from scalar flux and pressure, but only a single moment equation. 
Similarly, there is only a single boundary equation for both the normal component of the current and the boundary factor.
To solve the moment equations, closures for the pressure and boundary factor must be introduced.
The SMM uses additive closures given by
\begin{subequations}
    \begin{gather}
        \bP = \sum_{n=1}^{N} w_n (\bOmega_n \otimes \bOmega_n) \psi_n - \third \vec{I} \sum_{n=1}^{N} w_n \psi_n + \third \vec{I} \varphi = \bT(\psi) + \third \vec{I} \varphi \\
        B = \sum_{n=1}^{N} w_n \lvert \bOmega \cdot \bnormal \rvert \psi_n - f_b \sum_{n=1}^{N} w_n \psi_n +  f_b \varphi = \beta(\psi) + f_b \varphi,
    \end{gather}
\end{subequations}
where $\vec{I}$ is the identity tensor, the argument $\psi$ implies that the quantity is a function of the high-order transport solution, and
\begin{equation}
    f_b = \dfrac{\sum_{n=1}^{N} w_n \lvert \bOmega_n \cdot \bnormal \rvert}{\sum_{n=1}^{N} w_n} \approx \half,
\end{equation}
is used in lieu of 1/2 to account for the fact that $\lvert \bOmega \cdot \bnormal \rvert \psi$ cannot be integrated exactly by quadrature \cite{olivier2023smm,olivier2023vef}.
For the continuous equations, these closures are merely an algebraic manipulation since $\varphi - \sum_{n=1}^{N} w_n \psi_n = 0$.
When the transport and moment equations are inconsistently discretized, however, these closures introduce a difference on the order of the discretization error \cite{olivier2023smm}. 

After substituting the closures, the moment system becomes
\begin{subequations}
    \begin{equation}
        -\bdiv \dfrac{1}{3\sigma_{t,g}} \bgrad \varphi_g + \sigma_{r,g} \varphi_g = \sum_{g' \ne g} \sigma_{s,g'g} \varphi_{g'} + \dfrac{\chi_g}{k} \sum_{g'=1}^{G} \nu\sigma_{f,g'} \varphi_{g'} + \bdiv \dfrac{1}{\sigma_{t,g}} \bdiv \bT_g(\psi_g),
        \label{eq:smm}
    \end{equation}\vspace{-8ex}
    \begin{alignat}{2}
        &\bJ_g \cdot \bnormal = \beta_g(\psi_g) + f_b \varphi_g, \qquad &&\bx \in \partial\Gamma_\text{vac}, 
        \label{eq:smm-vac}\\
        &\bJ_g \cdot \bnormal = 0, \qquad &&\bx \in \partial\Gamma_\text{refl},
        \label{eq:smm-refl}
    \end{alignat}
    \label{eq:smm-full}
\end{subequations}
where $\bdiv \vec{I} \varphi = \bgrad \varphi$.
Note that this has the form of a standard multi-group diffusion $k$-eigenvalue problem with a volumetric and boundary source.
This moment system is not an approximation to the transport equation, but rather a reformulation of it. 
The moment system, therefore, still captures transport physics and is valid in non-diffusive regimes.
Further, when the transport solution is linearly anisotropic, the moment system reduces to exactly the diffusion equation with a Marshak boundary condition \cite{olivier2023smm}.
While outside the scope of this paper, the inclusion of anisotropic scattering would not materially alter this derivation.
Its inclusion would result in a source term in Eq. \eqref{eq:mg_first} which would be carried through to Eq. \eqref{eq:smm}.

\subsection{Second Moment Method Algorithm}

The transport $k$-eigenvalue problem can be written in operator notation as 
\begin{equation}
    L \bpsi = MSD \bpsi + \dfrac{1}{k} MFD \bpsi 
\end{equation}
where $\bpsi = [\psi_1, \ldots, \psi_G]$, and $\bphi = D \bpsi$ are the equivalent group-wise flux moments.
Here, $L$ is the transport operator, given by the left hand side of Eq. \eqref{eq:mgdo}, $M$ is the moment-to-discrete operator which maps a moment-based quantity to a direction-based quantity, $D$ is the discrete-to-moment operator, or the converse of $M$, $S$ is the scattering operator, and $F$ is the fission operator.

The SMM begins with a standard power iteration (PI) given by
\begin{equation}
    (L - MSD)\bpsi^{\ell+1/2} = \dfrac{1}{k^{\ell}} MFD \bpsi^{\ell},
\end{equation}
where $\ell$ is the power iteration index. 
This requires the solution to a fixed source problem driven by the previous power iteration fission source. 
This can be solved with standard techniques such as source iteration (SI), which lags the scattering source via
\begin{equation}
    L\bpsi^{\ell+1/2, i+1} = MSD\bpsi^{\ell+1/2, i} + \dfrac{1}{k^{\ell}} MFD \bpsi^{\ell},
\end{equation}
where $i$ is the source iteration index.
For the SMM, a single SI is performed to obtain
\begin{equation}
    \bpsi^{\ell+1/2} = L^{-1}\left( MSD + \dfrac{1}{k^\ell} MFD \right) \bpsi^{\ell},
\end{equation}
where $ \bpsi^{\ell+1/2} = \bpsi^{\ell+1/2,1}$, and $\bpsi^{\ell} = \bpsi^{\ell+1/2, 0}$.
With this, the SMM closures from Eqs. \eqref{eq:smm} can be computed, the diffusion correction source assembled into $\mathcal{R}(\bpsi^{\ell+1/2})$, and the scalar flux mapped to the moment system for the initial guess.
This mapping can be performed via $\bvarphi = \mathcal{P} D_0 \bpsi = \mathcal{P} \bphi_0$, where $D_i$ is the moment-to-discrete operator for only the $i$th moment and $\mathcal{P}$ is a conservative mapping from the transport discretization to the moment system.

Now, everything necessary to close the moment system is available and the moment system can be solved.
In operator notation, the moment system is given by
\begin{equation}
    \mathcal{D} \bvarphi = \mathcal{S} \bvarphi + \dfrac{1}{\lambda} \mathcal{F} \bvarphi + \mathcal{R}(\bpsi),
\end{equation}
where $\mathcal{D}$ is the diffusion operator, given by the left hand side of Eq. \eqref{eq:smm}, $\mathcal{S}$ and $\mathcal{F}$ are the diffusion scattering and fission operators, respectively, and $\lambda$ is the eigenvalue for the moment system.
The moment system, like the transport equation, can be solved using the power method via
\begin{equation}
    (\mathcal{D} - \mathcal{S}) \bvarphi^{m+1} = \dfrac{1}{\lambda^{m}} \mathcal{F} \bvarphi^{m} + \mathcal{R}(\bpsi^{\ell+1/2}),
\end{equation}
Because $\mathcal{D}$ is SPD while $\mathcal{D} - \mathcal{S}$ may not be, in order to use highly efficient linear solvers, the SI procedure is again employed to obtain
\begin{equation}
    \mathcal{D} \bvarphi^{m+1, j+1} = \mathcal{S} \bvarphi^{m+1, j} + \dfrac{1}{\lambda^{m}} \mathcal{F} \bvarphi^{m} + \mathcal{R}(\bpsi^{\ell+1/2}). 
\end{equation}
Once the inner iterations have converged or are terminated, the eigenvalue can be updated via
\begin{equation}
    \lambda^{m+1} = \dfrac{\lVert \mathcal{F} \bvarphi^{m+1} \rVert}{\lVert \mathcal{F} \bvarphi^{m} \rVert} \lambda^{m}.
\end{equation}
Once the power iterations have converged or are terminated, the transport eigenvalue takes the value of the diffusion eigenvalue, or $k^{\ell+1} = \lambda^{m+1}$, and the transport scalar flux becomes the diffusion scalar flux, or $\bphi_0^{\ell+1} = \widehat{\mathcal{P}} \bvarphi^{m+1}$, where $\widehat{\mathcal{P}}$ is the converse of $\mathcal{P}$.
The outer transport iterations converge once the the transport eigenvalue, and optionally the eigenvalue shape, meet the criteria
\begin{equation}
    \dfrac{\lvert k^{\ell+1} - k^\ell \rvert}{\lvert k^\ell \rvert} < \epsilon_k, \qquad \lVert \bphi_{0}^{\ell+1} - \bphi_{0}^{\ell} \rVert_{\ell_2} < \epsilon_{\bphi_0}.
\end{equation}
     
\section{Moment System Discretization}
\label{sec:discretization}

In this work, Eq. \eqref{eq:mgdo} is discretized with the PWLD finite element discretization of Bailey \cite{bailey2008piecewise}, and the moment system with a PWLC finite element discretization.
For a more detailed description of the PWL discretizations, the curious reader should consult \cite{bailey2008piecewise,vermaak2021massively} and the references therein.
For the purposes of this work, the discretization is considered in a general finite element framework.

For clarity, the discretization is presented for a single group and follows closely to Olivier et al. \cite{olivier2023vef}.
Restated, moment equations for a single group are given by
\begin{subequations}
    \begin{gather}
        \bdiv \bJ + \sigma_r \varphi = Q_0
        \label{eq:zeroth}\\
        \third \bgrad \varphi + \sigma_t \bJ = -\bdiv \bT,
        \label{eq:first} \\
        \bJ \cdot \bnormal = \beta + f_b \varphi, \qquad \bx \in \partial\Gamma_\text{vac} 
        \label{eq:moment-vac}\\
        \bJ \cdot \bnormal = 0, \qquad \bx \in \partial\Gamma_\text{refl}
        \label{eq:moment-refl}
    \end{gather}
\end{subequations}
where $Q_0$ contains the isotropic scattering and fission sources.
Here, the moment system unknowns $\varphi$ and $\bJ$ are continuous and the closure from transport $\bT$ is discontinuous.
While the current is trivially eliminated from the zeroth moment equation from the continuum moment system in Section \ref{sec:method}, the elimination of the current in the discrete system is complicated by the discontinuous transport discretization. 

Let the domain $\Gamma$ be subdivided into a set of arbitrary cells $\mathcal{T}$ such that $\Gamma = \bigcup_{K \in \mathcal{T}} K$.
Further let $\mathcal{F}$ be the collection of all unique faces, $\mathcal{F}_0 = \mathcal{F} \setminus \partial\Gamma$ be the set of unique interior faces,
and $\mathcal{F}_\text{vac} = \partial\Gamma_\text{vac}$ denote the set of unique boundary faces that belong to a vacuum boundary.
Multiplying Eq. \eqref{eq:zeroth} by continuous test function $u$ and integrating over the domain gives
\begin{equation}
    - \int_\mathcal{T} \bgrad u \cdot \bJ dV + f_b \int_{\mathcal{F}_\text{vac}} u \varphi  dS + \int_{\mathcal{F}_\text{vac}} u \beta dS + \int_\mathcal{T} \sigma_r u \varphi dV = \int_\mathcal{T} u Q_0 dV,
\end{equation}
where the first three terms arise from integration by parts and Eqs. \eqref{eq:moment-vac} and \eqref{eq:moment-refl}.
Because Eq. \eqref{eq:first} is a vector equation, it must be tested with a vector-valued test function $\vec{v}$.
Multiplying by $\vec{v}$ and integrating over cell $K$ gives
\begin{equation}
    \third \int_K \vec{v} \cdot \bgrad_K \varphi dV + \int_K \sigma_t \vec{v} \cdot \bJ dV = -\int_{\partial K} \vec{v} \cdot \widehat{\bT} \bnormal dS + \int_K \bgrad_K \vec{v} \colon \bT dV,
\end{equation}
where $\bgrad_K$ is the local gradient on cell $K$, $\widehat{\bT}$ is a numerical flux, and the tensor contraction operator is defined by
\begin{equation}
    \vec{A} \colon \vec{B} = \sum_{i=1}^{d} \sum_{j=1}^{d} A_{ij} B_{ij}.
\end{equation}
The integration is done only over cell $K$ in this case because particular care must be taken in handling surface discontinuities.
For continuous functions, these internal surface integrals cancel when summing over cells, leaving only boundary integrals.
For discontinuous functions, however, a numerical flux must be introduced to make the function single-valued at interior interfaces.
The particular definition of the numerical flux is arbitrary and left unspecified.
A second application of integration by parts yields
\begin{equation}
    \third \int_K \vec{v} \cdot \bgrad_K \varphi dV + \int_K \sigma_t \vec{v} \cdot \bJ dV = \int_{\partial K} \vec{v} \cdot \left( \bT \bnormal - \widehat{\bT} \bnormal \right) dS - \int_K \vec{v} \cdot \bgrad_K \cdot \bT dV,
\end{equation}
where a numerical flux is no longer required due to the integration by parts on the local gradient $\bgrad_K$.

In order to obtain a discrete elimination for the current, this must be extended to the entire domain $\Gamma$.
Summing over all cells $K$, this becomes
\begin{equation}
    \begin{aligned}
        \third \int_\mathcal{T} \vec{v} \cdot \bgrad &\varphi dV + \int_\mathcal{T} \sigma_t \vec{v} \cdot \bJ dV = \\ &\int_{\mathcal{F}} \jump{\vec{v}} \cdot \average{\bT \bnormal - \widehat{\bT} \bnormal} dS + \int_{\mathcal{F}_0} \average{\vec{v}} \cdot \jump{\bT \bnormal - \widehat{\bT} \bnormal} dS - \int_\mathcal{T} \vec{v} \cdot \bgrad_h \cdot \bT dV,
    \end{aligned}
\end{equation}
where $\bgrad_h$ is a generalization of $\bgrad_K$ to $\mathcal{T}$, and it can be shown that
\begin{equation}
    \sum_K \int_{\partial K} u \vec{v} \cdot \bnormal dS = \int_{\mathcal{F}} \jump{u} \average{\vec{v} \cdot \bnormal} dS + \int_{\mathcal{F}_0} \jump{\vec{v} \cdot \bnormal} \average{u} dS, 
\end{equation}
where
\begin{equation*}
    \jump{u} = u^+ - u^-, \qquad \average{u} = \half (u^+ + u^-),
\end{equation*}
are jump and average operators, respectively.
The $\pm$ notation denotes cell ownership at a discontinuity based on the orientation of the cell relative to a unique face normal.
Because $\vec{v}$ is continuous, the previous result can be reduced to
\begin{equation}
    \third \int_\mathcal{T} \vec{v} \cdot \bgrad \varphi dV + \int_\mathcal{T} \sigma_t \vec{v} \cdot \bJ dV = \int_{\mathcal{F}_0} \average{\vec{v}} \cdot \jump{\bT \bnormal - \widehat{\bT} \bnormal} dS - \int_\mathcal{T} \vec{v} \cdot \bgrad_h \cdot \bT dV.
\end{equation}
Choosing a conservative numerical flux such that $\jump{\widehat{\bT} \bnormal} = 0$, this becomes
\begin{equation}
    \third \int_\mathcal{T} \vec{v} \cdot \bgrad \varphi dV + \int_\mathcal{T} \sigma_t \vec{v} \cdot \bJ dV = \int_{\mathcal{F}_0} \average{\vec{v}} \cdot \jump{\bT \bnormal} dS - \int_\mathcal{T} \vec{v} \cdot \bgrad_h \cdot \bT dV.
\end{equation}

To obtain a discrete strong form of the current, all integrals must be over the domain $\mathcal{T}$ and have an integrand of the form $\vec{v} \cdot \vec{w}$ where $\vec{w}$ is an arbitrary vector-valued quantity.
This can be accomplished by defining an abstract lifting operator such that
\begin{equation}
    \int_\mathcal{T} \sigma_t \vec{v} \cdot \vec{r}\left( \jump{\bT \bnormal} \right) dV = \int_{\mathcal{F}_0} \average{\vec{v}} \cdot \jump{\bT \bnormal} dS,
    \label{eq:lifting}
\end{equation}
where $\vec{r}(\cdot)$ is non-zero only on $\mathcal{F}_0$.
Plugging this in gives
\begin{equation}
    \int_\mathcal{T} \sigma_t \vec{v} \cdot \bJ dV = -\third \int_\mathcal{T} \vec{v} \cdot \bgrad \varphi dV + \int_\mathcal{T} \sigma_t \vec{v} \cdot \vec{r}\left( \jump{\bT \bnormal} \right) dV - \int_\mathcal{T} \vec{v} \cdot \bgrad_h \cdot \bT dV,
\end{equation}
which results in the strong form current
\begin{equation}
    \bJ = \vec{r}\left( \jump{\bT \bnormal} \right) - \dfrac{1}{3\sigma_t} \bgrad \varphi - \dfrac{1}{\sigma_t} \bgrad_h \cdot \bT.
    \label{eq:strong_current}
\end{equation}
This result is equivalent to the standard definition of the current from the continuum moment equations within cells, but includes an addition term to account for discontinuities in $\bT$ on interior faces.
It should be emphasized that the lifting operator is simply a mathematical abstraction used to recast the surface integration as a volumetric integration, in turn allowing for the above definition.
Using the definitions of the lifting operator and the strong form current from Eqs. \eqref{eq:lifting} and \eqref{eq:strong_current}, the weak form of the moment system is given by
\begin{equation}
    \begin{aligned}
        \int_\mathcal{T} \bgrad u \cdot &\dfrac{1}{3\sigma_t} \bgrad \varphi dV + \int_\mathcal{T} \sigma_r u \varphi dV + f_b \int_{\mathcal{F}_\text{vac}} u \varphi dS = \\ &\int_\mathcal{T} u Q_0 dV - \int_{\mathcal{F}_\text{vac}} u \beta dS - \int_\mathcal{T} \bgrad u \cdot \dfrac{1}{\sigma_t} \bgrad_h \cdot \bT dV + \int_{\mathcal{F}_0} \average{\dfrac{\bgrad u}{\sigma_t}} \cdot \jump{\bT \bnormal} dS.
    \end{aligned}
\end{equation}
From this, the multi-group weak form can be trivially defined as
\begin{equation}
    \begin{aligned}
            &\int_\mathcal{T} \bgrad u_g \cdot \dfrac{1}{3\sigma_{t,g}} \bgrad \varphi_g dV + \int_\mathcal{T} \sigma_{r,g} u_g \varphi_g dV + f_b \int_{\mathcal{F}_\text{vac}} u_g \varphi_g dV = \\ &\qquad \int_\mathcal{T} \sum_{g' \ne g} \sigma_{s,g'g} u_g \varphi_{g'} + \int_\mathcal{T} \dfrac{\chi_g}{k} \sum_{g'=1}^{G} \nu\sigma_{f,g'} u_g \varphi_{g'} dV - \int_{\mathcal{F}_\text{vac}} u_g \beta_g dS \\ &\qquad\qquad - \int_\mathcal{T} \bgrad u \cdot \dfrac{1}{\sigma_{t,g}} \bgrad_h \cdot \bT_g dV + \int_{\mathcal{F}_0} \average{\dfrac{\bgrad u_g}{\sigma_{t,g}}} \cdot \jump{\bT_g \bnormal} dS, \qquad g = 1, \ldots, G,
    \end{aligned}
\end{equation}
where $\vec{u} = [u_1, \ldots, u_G]$ is a group-wise scalar test function.
 
\section{Results}
\label{sec:results}

In this section, the SMM is tested on two different problems. 
Its performance characterized relative to standard PI, a JFNK method \cite{knoll2011acceleration,park2012nonlinear}, and the Barbu-Adams method (DSA) \cite{barbu2019linear,barbu2023convergence}.
The JFNK method is a nonlinear formulation of the $k$-eigenvalue problem and the DSA discretization is algebraically consistent with the DO discretization.
Unless otherwise specified, all results use an outer iteration tolerance of $10^{-8}$.
For DSA and the SMM, the diffusion outer iterations use the same tolerance and a maximum of 50 iterations are allowed.

The SMM was implemented in the open-source massively parallel radiation transport code OpenSn (\href{https://github.com/Open-Sn/opensn}{https://github.com/Open-Sn/opensn}) \cite{vermaak2021massively}, which heavily utilizes the PETSc library \cite{petsc-user-ref}.
The PETSc GMRES solver is used for the inner transport iterations for both PI and JFNK methods, while Richardson iterations are used for the SMM and DSA \cite{petsc-user-ref}.
The JFNK method uses the PETSc scalable nonlinear equations solver (SNES) \cite{petsc-user-ref}. 
The low-order systems for DSA and the SMM are preconditioned with BoomerAMG from the \emph{hypre} library \cite{falgout2002hypre}, and solved with conjugate gradient.

\subsection{Simple 1D Example}

\begin{figure}[!htb]
    \centering
    \includegraphics[width=0.8\textwidth]{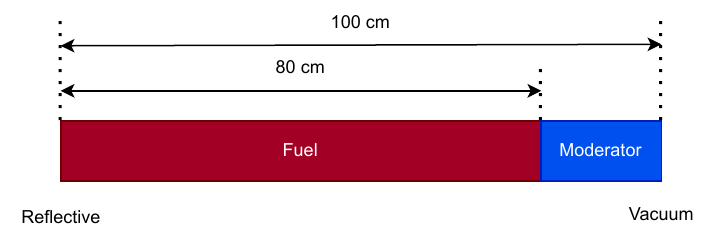}
    \caption{Geometry for the two region slab problem.}
    \label{fig:slab-geom}
\end{figure}

Consider a 100-cm slab where the first 80 cm from the left contain fissile material and the remainder is a moderating material, as depicted in Figure \ref{fig:slab-geom}.
The left boundary is defined as reflective and the right as vacuum.
The two-group cross sections for this problem are given in Table \ref{tab:slab-xs}. 
\begin{table}[!htb]
    \centering
    \caption{Cross section data for the two region slab problem.}
    \begin{tabular}{|c|c|c|c|c c| c |}
        \hline
        Material & Group & $\sigma_t$ & $\nu\sigma_f$ & $\sigma_{s,g1}$ & $\sigma_{s,g2}$ & c \\
        \hline
        \multirow{2}{*}{Fuel}       
        & 1 & 0.4241 & 0.02615 & 0.3944   & 0.0007568 & 0.932 \\
        & 2 & 0.7377 & 0.6285  & 0.001266 & 0.4021    & 0.547 \\
        \hline
        \multirow{2}{*}{Moderator}  \
        & 1 & 0.6823 & -       & 0.6531   & 0.0288    & 0.999 \\
        & 2 & 1.8690 & -       & 0.002143 & 1.857     & 0.995\\
        \hline
    \end{tabular}
    \label{tab:slab-xs}
\end{table}
While simple, this problem serves as a good basis to study the SMM.
Because of the optical thickness (high dominance ratio) and high scattering ratios, both the inner and outer iterations are slow to converge, making acceleration necessary.
For all results in this section, an $S_{64}$ angular quadrature is used.

\begin{figure}[!htb]
    \centering
    \includegraphics[width=0.8\textwidth]{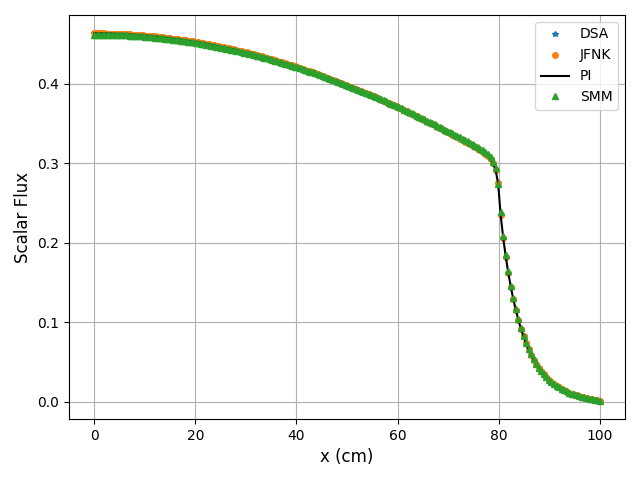}
    \caption{The eigenfunctions obtained from each method for the two region slab problem.}
    \label{fig:slab-shapes}
\end{figure}

In Table \ref{tab:slab-comp}, the iterative performance of each method is compared.
\begin{table}[!htb]
    \centering
    \caption{Two region slab problem performance in serial.}
    \begin{tabular}{|c|c|c|c|c|c|}
        \hline
        \# of Cells & Method & $k_\text{eff}$ & \# of Outers & \# of Sweeps \\
        \hline
        \multirow{5}{*}{500 Cells} 
        & PI    & 0.925447  &  189 & 4347 \\
        & JFNK  & 0.925447  &  5   & 198 \\
        & DSA   & 0.925447  &  14  & 14  \\
        & SMM   & 0.925461  &  13  & 13   \\
        \hline
        \multirow{5}{*}{5000 Cells}
        & PI    & 0.925447  &  189 & 4347 \\
        & JFNK  & 0.925447  &  5   & 213 \\
        & DSA   & 0.925447  &  13  & 13  \\
        & SMM   & 0.925447  &  12  & 12  \\
        \hline
    \end{tabular}
    \label{tab:slab-comp}
\end{table}
The SMM requires a comparable number of sweeps to DSA and dramatically fewer than both PI and JFNK.
Further, under mesh refinement, the number of sweeps remains roughly constant.
Because JFNK solves the exact transport $k$-eigenvalue problem and DSA is consistent with the transport discretization, each of these methods return the same eigenvalue as unaccelerated  transport (PI).
The algebraic inconsistency of the SMM discretization, on the other hand, results in a deviation.
As the mesh is refined, however, the SMM eigenvalue converges to the transport eigenvalue.
For the 500-cell case, there is a 14-pcm difference, whereas for the 5000-cell case, the SMM is within one pcm.
The eigenfunctions obtained from each method are plotted against one another in Figure \ref{fig:slab-shapes}, where all eigenfunctions are visually identical.

\begin{figure}[!htb]
    \centering
    \includegraphics[width=0.8\textwidth]{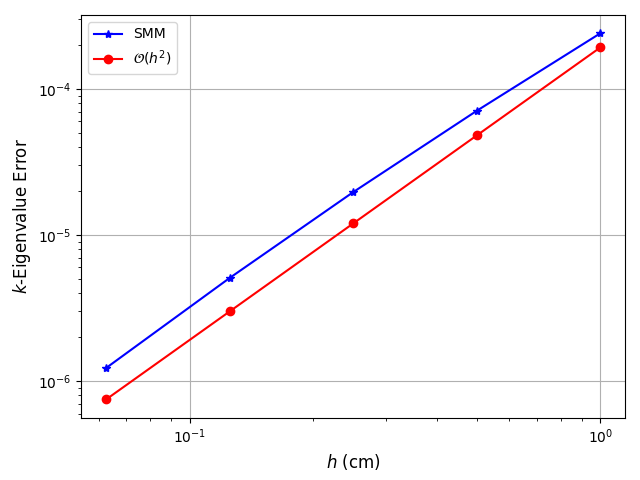}
    \caption{Error as a function of mesh size for the two region slab problem.}
    \label{fig:smm-convergence}
\end{figure}

The order of convergence for the SMM is estimated using 
\begin{equation}
    \mathcal{O} = \log_2\left( \dfrac{\lvert k_N - k_\text{ref} \rvert}{\lvert k_{2N} - k_\text{ref} \rvert} \right),
\end{equation}
where $k_\text{ref}$ is a reference eigenvalue on a fine mesh and $N$ is the number of cells on a coarse mesh.
Because both the transport and moment systems are discretized with second order methods, halving the mesh size should result in a reduction in error by a factor of four.
This behavior is observed in Table \ref{tab:slab-convergence}, where $e$ and $\mathcal{O}$ are the error and approximate order of convergence relative to the fine mesh SMM eigenvalue, respectively, and $e_\text{DO}$ and $\mathcal{O}_\text{DO}$ are defined equivalently for the fine mesh transport eigenvalue.
The reference PI and SMM eigenvalues, $0.92544736$ and $0.92544783$, respectively, were computed on a grid with 3200 cells.
Figure \ref{fig:smm-convergence} shows the eigenvalue error as a function of the cell width plotted alongside a line depicting true second order convergence.
\begin{table}[!htb]
    \centering
    \caption{Convergence study for the two region slab problem.}
    \begin{tabular}{|c|c|c|c|c|c|}
        \hline
        Cells & $k_\text{eff}$ & $e$ & $e_\text{DO}$ & $\mathcal{O}$ & $\mathcal{O}_\text{DO}$ \\
        \hline
        100  & 0.9256877 & 2.3980e-04 & 2.4008e-04 & -     & -  \\ 
        200  & 0.9255184 & 7.0527e-05 & 7.0807e-05 & 1.766 & 1.762 \\ 
        400  & 0.9254672 & 1.9385e-05 & 1.9666e-05 & 1.863 & 1.848 \\ 
        800  & 0.9254527 & 4.8051e-06 & 5.0855e-06 & 2.012 & 1.951 \\ 
        1600 & 0.9254488 & 9.5079e-07 & 1.2312e-06 & 2.337 & 2.046 \\ 
        \hline
    \end{tabular}
    \label{tab:slab-convergence}
\end{table}

\subsection{C5G7 Benchmark Problem}

\begin{figure}[!htb]
    \centering
    \includegraphics[width=0.8\linewidth]{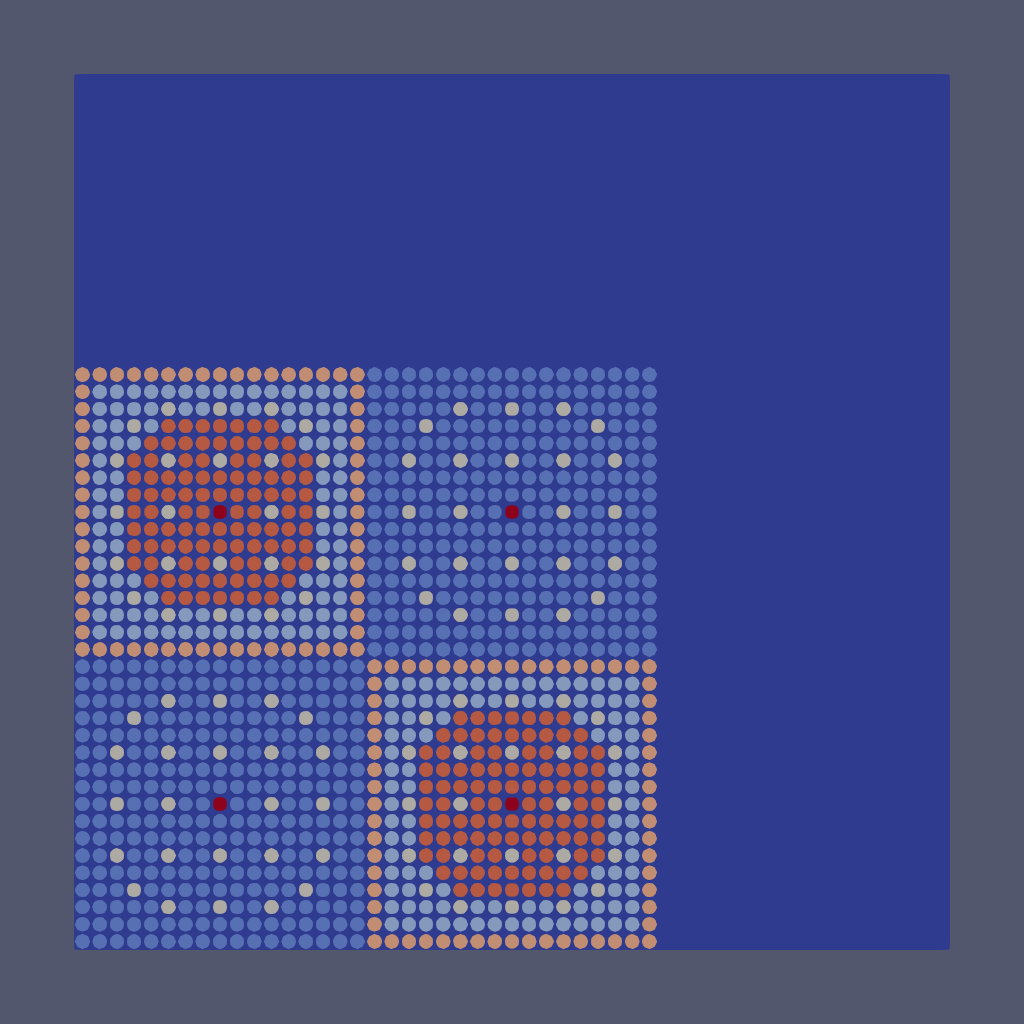}
    \caption{The geometry for the C5G7 benchmark problem. Each assembly has a central fission chamber (dark red) and a grid of guide tubes (tan). The bottom left and top right fuel assemblies are UO2 fuel and the top left and bottom right contain MOX fuels of various enrichment levels. From interior to exterior, the MOX enrichment levels are 8.7\% (red), 7.0\% (light blue), and 4.3\% (orange).}
    \label{fig:c5g7-geom}
\end{figure}

The C5G7 benchmark is a common problem used to characterize the performance of neutron transport methods.
This problem is a two-dimensional quarter-core $64.26 \times 64.26$-cm geometry with a two-by-two grid of $21.42 \times 21.42$-cm fuel assemblies surrounded by water, as shown in Figure \ref{fig:c5g7-geom}.
Each fuel assembly is a $17 \times 17$ lattice with a $1.26$-cm pitch and $0.54$-cm pin radius.
The bottom and left boundary conditions are reflective, while the top and right are vacuum.
The problem uses 7 energy groups and all scattering is isotropic.
More details regarding the geometry and cross sections can be found in \cite{lewis2001c5g7}.
The reference $k$-eigenvalue from a multi-group Monte-Carlo calculation is 1.186550 \cite{lewis2001c5g7}.

\begin{figure}[!htb]
    \centering
    \begin{subfigure}{0.45\linewidth}
        \centering
        \includegraphics[width=\linewidth]{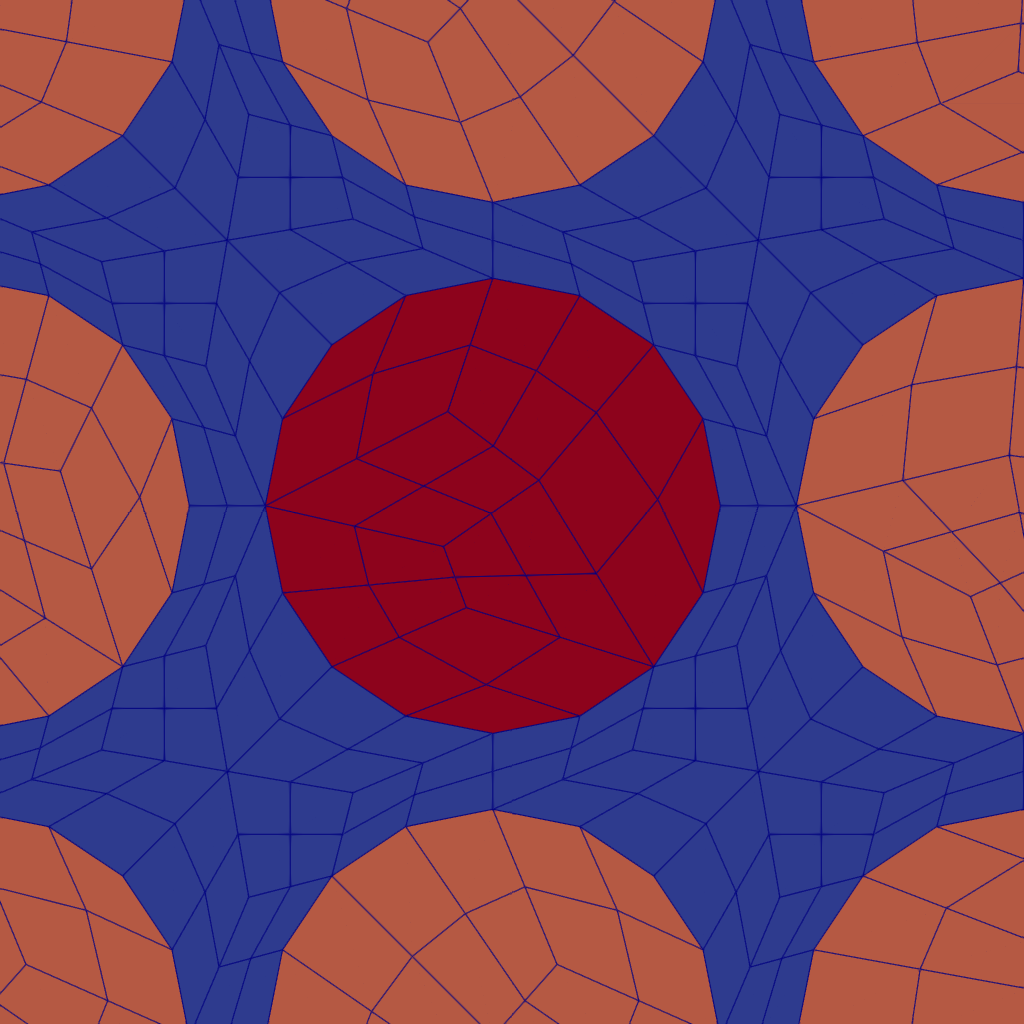}
     \end{subfigure}
     \begin{subfigure}{0.45\linewidth}
        \centering
        \includegraphics[width=\linewidth]{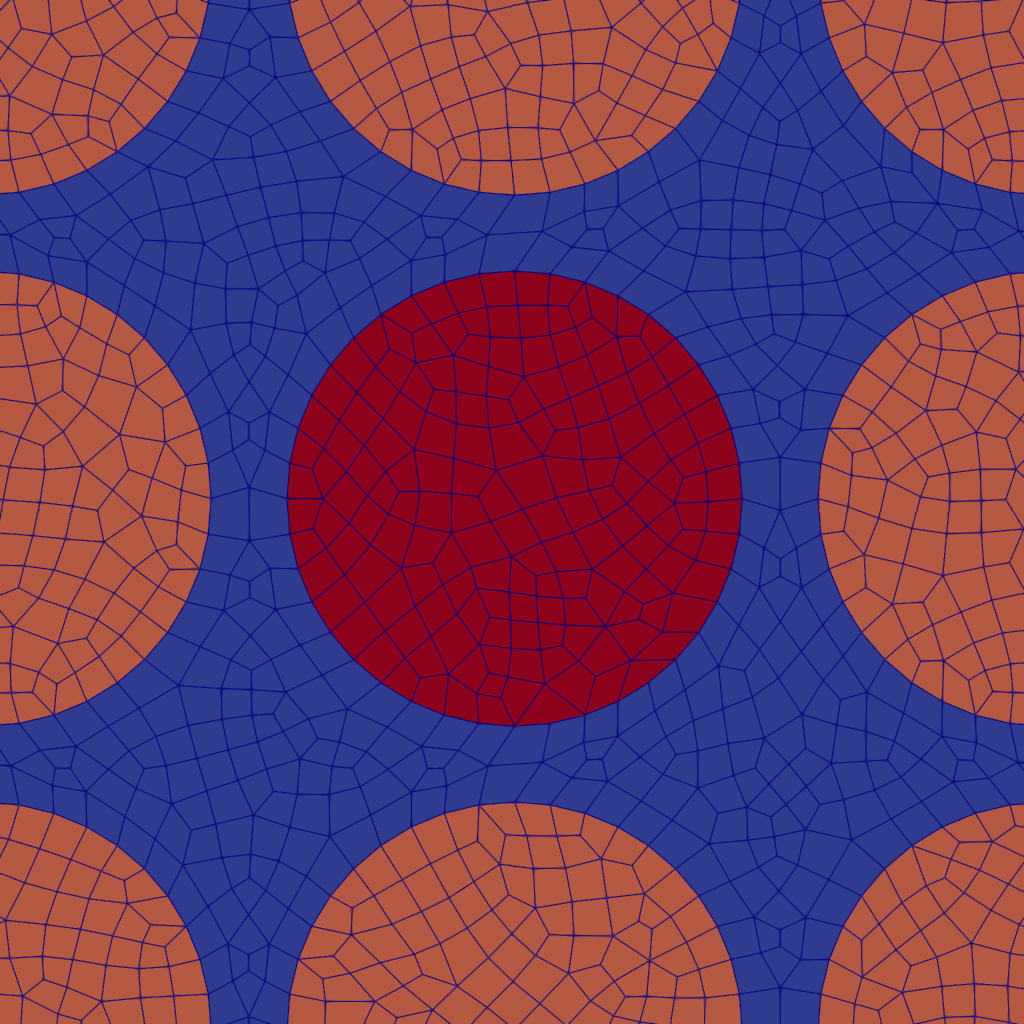}
     \end{subfigure}
     \caption{The coarse (left) and fine (right) C5G7 meshes zoomed in on a MOX fuel assembly.}
     \label{fig:c5g7-meshes}
\end{figure}

\begin{figure}
    \begin{subfigure}{0.45\linewidth}
        \centering
        \includegraphics[width=\linewidth]{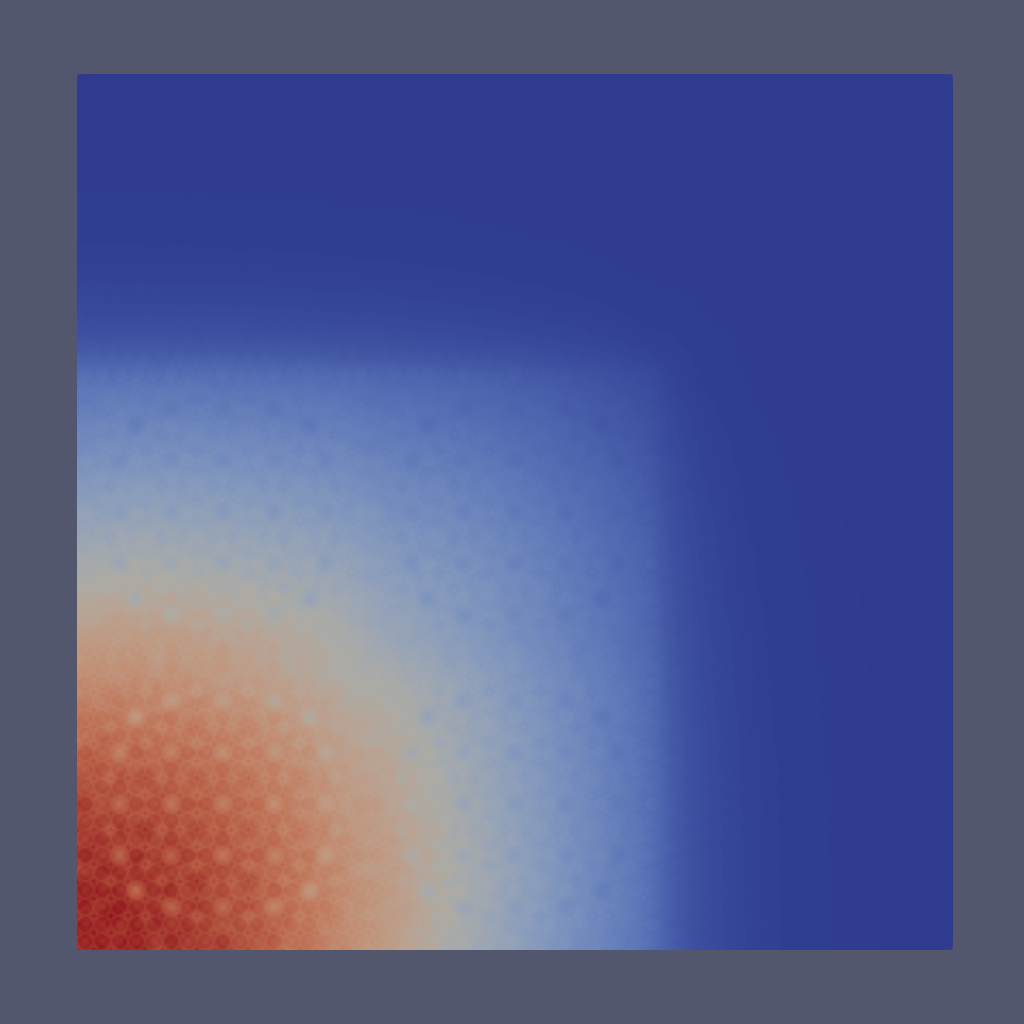}
    \end{subfigure}
    \begin{subfigure}{0.45\linewidth}
        \centering
        \includegraphics[width=\linewidth]{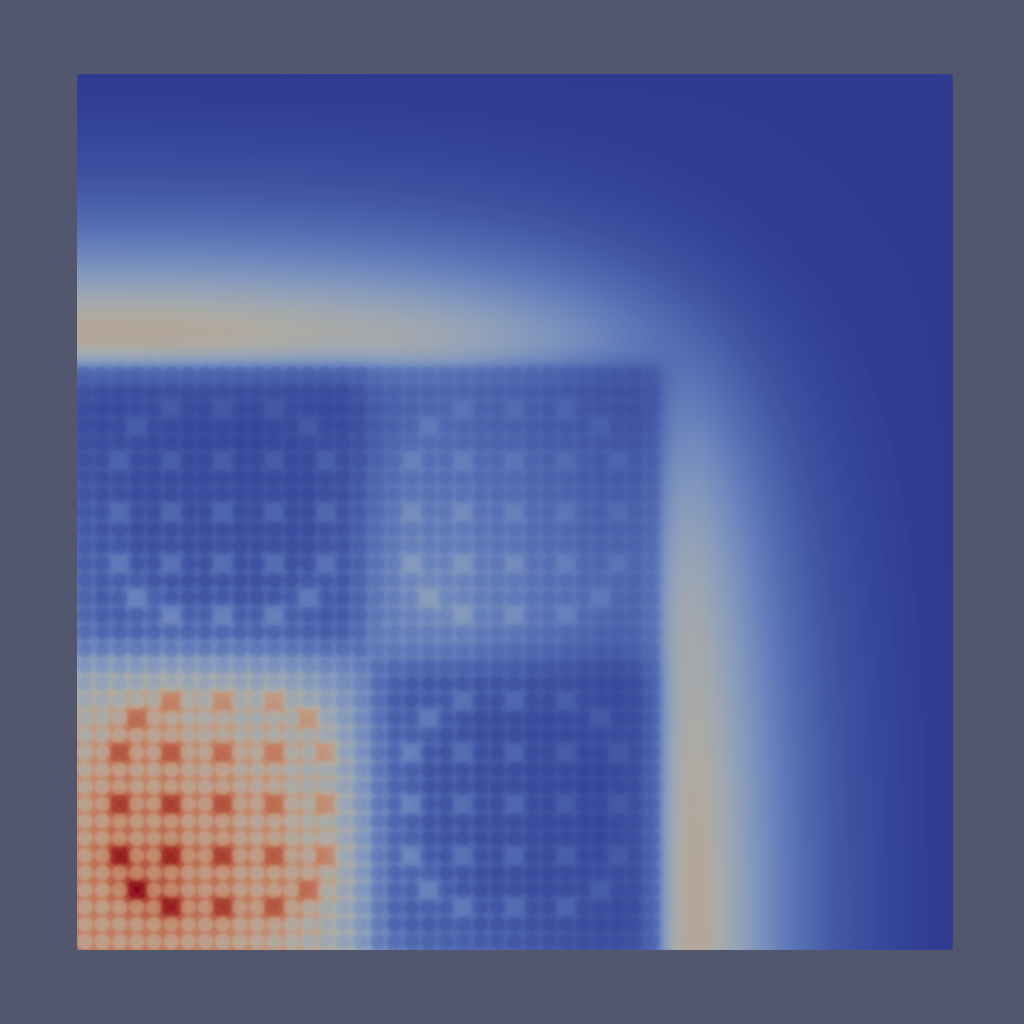}
    \end{subfigure}
    \caption{The scalar flux for the first (left) and last (right) group for the C5G7 benchmark problem.}
    \label{fig:c5g7-phi}
\end{figure}

Using \emph{gmsh} \cite{geuzaine2009gmsh}, the C5G7 benchmark problem is meshed using both unstructured quadrilaterals and triangles.
Two different meshes of varying refinement levels are used in this section.
The coarse mesh contains 83,264 cells, while the fine mesh contains 454,491.
With 7 energy groups and 16 directions, the total number of angular unknowns is then 37,287,040 and 203,611,968, respectively.
A zoomed in depiction of the meshes are given in Figure \ref{fig:c5g7-meshes}.
For reference, the first and last group scalar flux of the fundamental mode is given in Figure \ref{fig:c5g7-phi}.

Because of the difference introduced by the inconsistent discretization used in the SMM, a direct run time comparison between the different methods is inadequate.
Consider the figure of merit for efficiency,
\begin{equation}
    \eta = \dfrac{1}{|k - k_\text{ref}| T P},
\end{equation}
where $k$ is the simulation eigenvalue, $k_\text{ref}$ is the reference Monte-Carlo eigenvalue, $T$ is the clock time, and $P$ is the number of processors.
The purpose of this metric is to quantify the accuracy of a method per unit computational cost, where a larger value implies greater efficiency.
For two methods to have the same efficiency, the proportional difference in compute time must be the inverse of the proportional difference in error.
In other words, if one method has half the compute time and twice the error as another method, the methods have the same efficiency.
If two methods with the same efficiency were run for the same amount of time, two methods should obtain the same level of accuracy.

First, the performance of the SMM is compared to JFNK and DSA as a function of the number of directions used.
Table \ref{tab:c5g7-coarse} shows these results for the coarse mesh.
\begin{table}[!htb]
	\centering
	\caption{Coarse C5G7 benchmark problem performance with 16 MPI processes.}
	\begin{tabular}{|c|c|c|c|c|c|c|c|}
		\hline
		\# of Directions & Method & $k$  & \# of Outers & \# of Sweeps & Clock Time (s) & $\eta$\\
        \hline
		\multirow{3}{*}{4} 
        & JFNK  & 1.193239 & 7  & 133 & 123   & 0.076 \\
		& DSA   & 1.193239 & 16 & 16  & 1387  & 0.007 \\
		& SMM   & 1.193649 & 15 & 15  & 436   & 0.020 \\\hline
		\multirow{3}{*}{16} 
        & JFNK  & 1.192560 & 7  & 133 & 272   & 0.038\\
		& DSA   & 1.192559 & 13 & 13  & 1439  & 0.007\\
		& SMM   & 1.192709 & 14 & 14  & 366   & 0.028\\ \hline
		\multirow{3}{*}{64}
        & JFNK  & 1.192661 & 7  & 133 & 729  & 0.014 \\
		& DSA   & 1.192660 & 15 & 15  & 1745 & 0.006 \\
		& SMM   & 1.193023 & 14 & 14  & 423  & 0.023 \\\hline
        \multirow{3}{*}{256}
        & JFNK  & 1.193825 & 7  & 133 & 2305 & 0.004 \\
        & DSA   & 1.193825 & 14 & 14  & 1723 & 0.005 \\
        & SMM   & 1.194113 & 14 & 14  & 718  & 0.012\\\hline
	\end{tabular}
    \label{tab:c5g7-coarse}
\end{table}
The SMM universally converges with fewer sweeps than JFNK and with a comparable number to DSA.
Its compute time is also considerably faster than DSA in all cases, and JFNK when a realistic number of directions is used.
Because JFNK solves the transport $k$-eigenvalue problem directly, all nonlinear inner iterations require a sweep, making its computational cost scale with the number of directions.
For this reason, when a small number of directions are used and sweeps are relatively inexpensive, JFNK outperforms.
The SMM's significantly faster compute time compared to DSA is a direct result of having 3-4 times fewer interior unknowns due to the use of a continuous discretization.

The SMM eigenvalues are between 16-43 pcm greater than those of DSA and JFNK as a result of the inconsistent discretization.
Despite this, the figures of merit in Table \ref{tab:c5g7-coarse} show that the SMM is more efficient than DSA in all cases and JFNK for a realistic number of directions.
In other words, the gain in performance by using the SMM was greater than the accompanying loss in accuracy.

To further illustrate the concept of efficiency, consider the SMM and DSA results for the 64-direction problem.
The DSA and SMM errors are $e_\text{DSA} = 611$ pcm and $e_\text{SMM} = 1.06 e_\text{DSA} = 647$ pcm, respectively.
The compute times are $T_\text{DSA} = 1745$ s and $T_\text{SMM} = 0.24 T_\text{DSA} = 423$ s.
In this scenario, DSA required four times the computational cost to achieve a 6\% reduction in error relative to the SMM.
Instead of committing resources to run DSA on a coarse grid, if the SMM were run on a refined grid such that $T'_\text{SMM} = T_\text{DSA}$, a significant reduction in error would be expected because the method achieves greater accuracy per unit time.

As a larger test, the C5G7 benchmark problem is run on the fine mesh with 256 directions.
In this configuration, there are 3,257,791,488 unknowns.
The results of a strong scaling study shown in Table \ref{tab:c5g7-fine} were generated on the \emph{sawtooth} machine at Idaho National Laboratory (INL), which uses dual Xeon Platinum 8268 processors with 24 cores and 196 GB of memory.
\begin{table}[!htb]
    \centering
    \caption{Results for the refined C5G7 benchmark problem performance with 256 directions.}
    \begin{tabular}{|c|c|c|c|c|c|c|}
    \hline
        \# of Processes & Method & $k$ & \# of Outers & \# of Sweeps & Clock Time (s) & $\eta$ \\
        \hline
        \multirow{3}{*}{48}
        & JFNK  & 1.1873848 & 7  & 155 & 2889 & 0.009\\
        & DSA   & 1.1873847 & 12 & 12  & 1032 & 0.024\\
        & SMM    & 1.1874659 & 11 & 11  & 582  & 0.039\\
        \hline
        \multirow{3}{*}{96}
        & JFNK  & 1.1873848 & 7  & 155 & 1392 & 0.009\\
        & DSA   & 1.1873847 & 12 & 12  & 500  & 0.025\\
        & SMM    & 1.1874659 & 11 & 11  & 290  & 0.039\\
        \hline
        \multirow{3}{*}{192}
        & JFNK  & 1.1873848 & 7  & 155 & 851 & 0.007 \\
        & DSA   & 1.1873847 & 12 & 12  & 263 & 0.024\\
        & SMM    & 1.1874659 & 11 & 11  & 164 & 0.035 \\
        \hline
    \end{tabular}
    \label{tab:c5g7-fine}
\end{table}
Here, it is evident that the same relative performance between the three methods as before holds.
Further, the difference in the SMM eigenvalue to those of JFNK and DSA decreased to 8 pcm from 29 pcm on the coarse mesh.
This behavior is demonstrative of the behavior of the SMM under mesh refinement.
Because each method converges towards the true solution with the same order of accuracy, the difference between the SMM and other methods necessarily diminishes under mesh refinement at the same rate.

\begin{figure}
    \centering
    \includegraphics[scale=0.7]{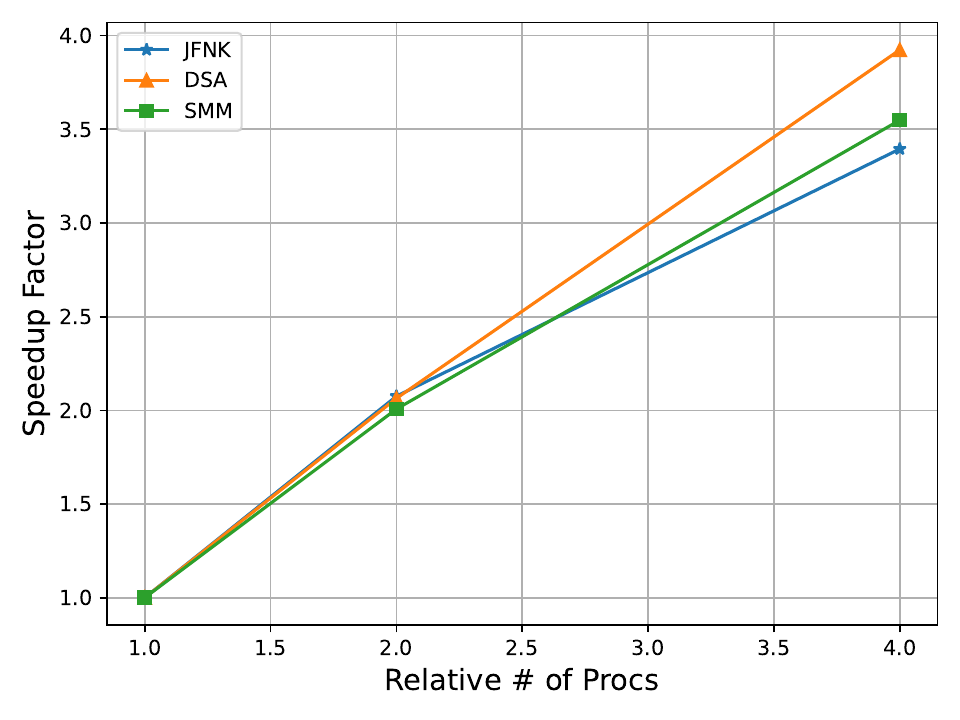}
    \caption{The strong scaling speedup factor as a function of the relative number of processors used for the fine mesh C5G7 benchmark problem with 256 directions.}
    \label{fig:strong-scaling}
\end{figure}

A plot of the strong scaling speedup factor is shown in Figure \ref{fig:strong-scaling}.
The strong scaling speedup factor is defined by
\begin{equation}
    \epsilon(n, n_b) = \dfrac{t(n_b)}{t(n)},
\end{equation}
where $n_b$ is the baseline number of processors, $n > n_b$ is the number of processors for a given simulation, and $t$ is the simulation time.
The ideal speedup factor is then $\epsilon(n, n_b) = n/n_b$.
All three methods scale quite well, with DSA very near the optimal of 4.0 and the SMM just above 3.5.
The slight degradation in parallel efficiency is primarily due to the existence of shared nodes in the continuous diffusion discretization.
This results in addition communications when projecting the transport solution to the diffusion discretization and during the linear diffusion solves.
The performance gains observed by using the SMM, however, dramatically outweigh this slight degradation.

In each of the above examples, the vast majority of the computational work of the SMM lies in the inner diffusion solves.
This is precisely the intention of methods like the SMM, QD, and DSA.
Scattering and fission are often slow to converge relative to the angular distribution.
By converging these physics in the low-order diffusion system, iterations progress more rapidly and the overall cost of the method becomes less sensitive to the number of directions.
For example, Table \ref{tab:c5g7-coarse} shows that there are very small differences in compute time between 4 and 64 directions.
In contrast, the 64-direction JFNK solution took roughly six times longer than the 4-direction solution.

\section{Conclusion}
\label{sec:conclusion}

This paper expanded the SMM $k$-eigenvalue acceleration scheme introduced in Tutt et al. \cite{tutt2023second} and Woodsford et al. \cite{woodsford2023variant} to general geometry, unstructured meshes, more advanced spatial discretizations, and multiple energy groups.
While it was previously thought that linear acceleration schemes lead to ill-posed eigenvalue problems, Barbu and Adams \cite{barbu2023convergence} empirically showed otherwise.
These results, like those presented in Tutt et al. \cite{tutt2023second} and Woodsford et al. \cite{woodsford2023variant}, help support that conclusion.
Unlike QD, the SMM results in an SPD system, and unlike DSA, a consistent diffusion discretization is not required for stability.
This allows for not only the use of highly efficient linear solvers, but for simpler diffusion discretizations with significantly fewer unknowns.
Each of these characteristics resulted in the SMM achieving significantly greater computational efficiency than its alternatives.

For a simple multi-material, multi-group slab problem with a high dominance ratio, the SMM required a negligible number of sweeps compared to both PI and JFNK, and was comparable to DSA.
While a small deviation in the eigenvalue is introduced due to the inconsistent SMM discretization, the deviation vanished as the mesh was refined, as expected.
This was further evidenced by the second order convergence of the SMM eigenvalue with respect to both fine mesh SMM and transport eigenvalues.

With the effectiveness of the method on a simple one-dimensional multi-group problem established, the method was then tested on the C5G7 benchmark problem \cite{lewis2001c5g7} for a more robust performance characterization.
Compared to DSA, the SMM had significantly lower compute times due to the reduced number of diffusion unknowns as a result of using an inconsistent discretization.
Provided that a realistic number of directions were used, the SMM was also dramatically faster than JFNK due to the small number of sweeps required.
Using a figure of merit for computational efficiency, the SMM sustained a factor of roughly 1.6 and 4 times greater efficiency than DSA and JFNK, respectively, on a problem with $3.25 \times 10^9$ unknowns solved on up to 192 processes.
Lastly, under mesh refinement, the deviation in the SMM eigenvalue from those of DSA and JFNK diminished considerably, as expected.

There are many potential extensions of this work.
First, the method could be extended to include anisotropic scattering in a relatively straightforward manner.
Further, a comparison to other nonlinear acceleration schemes such as CMFD or NDA would provide additional insights into the relative performance of the method.
Whereas the JFNK method uses a nonlinear formulation of the transport eigenvalue problem, a similar approach to that of Prince et al. \cite{prince2020diffusion} could be adopted where the diffusion eigenvalue problem is cast as a nonlinear problem and solved with a JFNK solver.
As discussed in Tutt et al. \cite{tutt2023second} and Woodsford et al. \cite{woodsford2023variant}, developing a coarse mesh version of the SMM would yield many advantageous characteristics.
Not only would a coarse mesh SMM yield an SPD system, but it would even further reduce the size of the linear system solved in the inner-most iterations of the algorithm.

\pagebreak
\section*{Acknowledgments}

This work was supported by the U.S. Department of Energy through the Los Alamos National Laboratory. Los Alamos National Laboratory is operated by Triad National Security, LLC, for the National Nuclear Security Administration of U.S. Department of Energy (Contract No. 89233218CNA000001).

This research made use of Idaho National Laboratory’s High Performance Computing systems located at the Collaborative Computing Center and supported by the Office of Nuclear Energy of the U.S. Department of Energy and the Nuclear Science User Facilities under Contract No. DE-AC07-05ID14517.

\pagebreak
\bibliographystyle{ans_js}                                                                           %custom ANS journal submission template bibliography style
\bibliography{SMM}

\begin{thebibliography}{10}
\newcommand{\enquote}[1]{``#1''}
\providecommand{\url}[1]{\texttt{#1}}
\providecommand{\urlprefix}{URL }
\expandafter\ifx\csname urlstyle\endcsname\relax
  \providecommand{\doi}[1]{doi:\discretionary{}{}{}#1}\else
  \providecommand{\doi}{doi:\discretionary{}{}{}\begingroup
  \urlstyle{rm}\Url}\fi

\bibitem{lewis1976comparison}
\textsc{E.~Lewis} and \textsc{W.~Miller~Jr}, \enquote{Comparison of P1
  synthetic acceleration techniques,} \emph{Transactions of the American
  Nuclear Society}, \textbf{23} (1976).

\bibitem{adams2002fast}
\textsc{M.~L. Adams} and \textsc{E.~W. Larsen}, \enquote{Fast iterative methods
  for discrete-ordinates particle transport calculations,} \emph{Progress in
  nuclear energy}, \textbf{40}, \emph{1}, 3 (2002).

\bibitem{chacon2017multiscale}
\textsc{L.~Chacon}, \textsc{G.~Chen}, \textsc{D.~A. Knoll}, \textsc{C.~Newman},
  \textsc{H.~Park}, \textsc{W.~Taitano}, \textsc{J.~A. Willert}, and
  \textsc{G.~Womeldorff}, \enquote{Multiscale high-order/low-order (HOLO)
  algorithms and applications,} \emph{Journal of Computational Physics},
  \textbf{330}, 21 (2017).

\bibitem{goldin1964quasi}
\textsc{V.~Y. Gol'din}, \enquote{A quasi-diffusion method of solving the
  kinetic equation,} \emph{Zhurnal Vychislitel'noi Matematiki i Matematicheskoi
  Fiziki}, \textbf{4}, \emph{6}, 1078 (1964).

\bibitem{olivier2023smm}
\textsc{S.~Olivier} and \textsc{T.~S. Haut}, \enquote{High-Order Finite Element
  Second Moment Methods for Linear Transport,} \emph{arXiv preprint
  arXiv:2304.07386} (2023).

\bibitem{tutt2023second}
\textsc{J.~Tutt}, \textsc{C.~Woodsford}, and \textsc{J.~E. Morel}, \enquote{A
  Second-Moment Method for $k$-Eigenvalue Acceleration,} \emph{Proceedings of
  M\&C 2023 International Conference on Mathematics and Computational Methods
  Applied to Nuclear Science and Engineering} (2023).

\bibitem{woodsford2023variant}
\textsc{C.~Woodsford}, \textsc{J.~Tutt}, and \textsc{J.~E. Morel}, \enquote{A
  Variant of the Second-Moment Method for $k$ Eigenvalue Calculations,}
  Submitted for publication in Nuclear Science and Engineering (2023).

\bibitem{ruge1987algebraic}
\textsc{J.~W. Ruge} and \textsc{K.~St{\"u}ben}, \enquote{Algebraic multigrid,}
  \emph{Multigrid methods}, 73--130, SIAM.

\bibitem{reed1971effectiveness}
\textsc{W.~H. Reed}, \enquote{The effectiveness of acceleration techniques for
  iterative methods in transport theory,} \emph{Nuclear Science and
  Engineering}, \textbf{45}, \emph{3}, 245 (1971).

\bibitem{alcouffe1977diffusion}
\textsc{R.~E. Alcouffe}, \enquote{Diffusion synthetic acceleration methods for
  the diamond-differenced discrete-ordinates equations,} \emph{Nuclear Science
  and Engineering}, \textbf{64}, \emph{2}, 344 (1977).

\bibitem{warsa2002fully}
\textsc{J.~S. Warsa}, \textsc{T.~A. Wareing}, and \textsc{J.~E. Morel},
  \enquote{Fully consistent diffusion synthetic acceleration of linear
  discontinuous SN transport discretizations on unstructured tetrahedral
  meshes,} \emph{Nuclear science and engineering}, \textbf{141}, \emph{3}, 236
  (2002).

\bibitem{adams1992diffusion}
\textsc{M.~L. Adams} and \textsc{W.~R. Martin}, \enquote{Diffusion synthetic
  acceleration of discontinuous finite element transport iterations,}
  \emph{Nuclear science and engineering}, \textbf{111}, \emph{2}, 145 (1992).

\bibitem{wang2010diffusion}
\textsc{Y.~Wang} and \textsc{J.~C. Ragusa}, \enquote{Diffusion synthetic
  acceleration for high-order discontinuous finite element SN transport schemes
  and application to locally refined unstructured meshes,} \emph{Nuclear
  science and engineering}, \textbf{166}, \emph{2}, 145 (2010).

\bibitem{olivier2023vef}
\textsc{S.~Olivier}, \textsc{W.~Pazner}, \textsc{T.~S. Haut}, and \textsc{B.~C.
  Yee}, \enquote{A family of independent Variable Eddington Factor methods with
  efficient preconditioned iterative solvers,} \emph{Journal of Computational
  Physics}, \textbf{473}, 111747 (2023).

\bibitem{cefus1989stability}
\textsc{G.~R. Cefus} and \textsc{E.~W. Larsen}, \enquote{Stability analysis of
  the quasideffusion and second moment methods for iteratively solving
  discrete-ordinates problems,} \emph{Transport Theory and Statistical
  Physics}, \textbf{18}, \emph{5-6}, 493 (1989).

\bibitem{smith2002full}
\textsc{K.~S. Smith}, \enquote{Full-core, 2-D, LWR core calculation with
  CASMO-4E,} \emph{Proc. Int. Conf. on the New Frontiers of Nuclear Technology:
  Reactor Physics, Safety and High-Performance Computing (PHYSOR2002), Oct.
  7-10, 2002, Seoul, Korea} (2002).

\bibitem{willert2014comparison}
\textsc{J.~Willert}, \textsc{H.~Park}, and \textsc{D.~A. Knoll}, \enquote{A
  comparison of acceleration methods for solving the neutron transport
  k-eigenvalue problem,} \emph{Journal of Computational Physics}, \textbf{274},
  681 (2014).

\bibitem{ramm2001simple}
\textsc{A.~G. Ramm}, \enquote{A simple proof of the Fredholm alternative and a
  characterization of the Fredholm operators,} \emph{The American Mathematical
  Monthly}, \textbf{108}, \emph{9}, 855 (2001).

\bibitem{gelbard1982acceleration}
\textsc{E.~Gelbard}, \textsc{C.~Adams}, \textsc{E.~Larsen}, and
  \textsc{D.~McCoy}, \enquote{Acceleration of transport eigenvalue problems,}
  \emph{Transactions of the American Nuclear Society}, \textbf{41} (1982).

\bibitem{dodson2019sda}
\textsc{Z.~Dodson}, \textsc{N.~Adamowicz}, \textsc{B.~Kochunas}, and
  \textsc{E.~Larsen}, \enquote{Sda: A semilinear cmfd-like transport
  acceleration scheme without d,} \emph{International Conference on Mathematics
  and Computational Methods Applied to Nuclear Science and Engineering (M\&C
  2019)}, 1714--1723 (2019).

\bibitem{barbu2019linear}
\textsc{A.~Barbu} and \textsc{M.~Adams}, \enquote{A linear
  diffusion-acceleration method for k-eigenvalue transport,} \emph{Proceedings
  of International Conference on Mathematics and Computational Methods Applied
  to Nuclear Science and Engineering, M\&C} (2019).

\bibitem{barbu2023convergence}
\textsc{A.~P. Barbu} and \textsc{M.~L. Adams}, \enquote{Convergence Properties
  of a Linear Diffusion-Acceleration Method for k-Eigenvalue Transport
  Problems,} \emph{Nuclear Science and Engineering}, \textbf{197}, \emph{4},
  517 (2023).

\bibitem{bailey2008piecewise}
\textsc{T.~S. Bailey}, \emph{The piecewise linear discontinuous finite element
  method applied to the RZ and XYZ transport equations}, Texas A\&M University
  (2008).

\bibitem{lewis2001c5g7}
\textsc{E.~Lewis}, \textsc{M.~Smith}, \textsc{N.~Tsoulfanidis},
  \textsc{G.~Palmiotti}, \textsc{T.~Taiwo}, and \textsc{R.~Blomquist},
  \enquote{Benchmark specification for Deterministic 2-D/3-D MOX fuel assembly
  transport calculations without spatial homogenization (C5G7 MOX),}
  \emph{NEA/NSC}, \textbf{280}, 2001 (2001).

\bibitem{lewis1984computational}
\textsc{E.~E. Lewis} and \textsc{W.~F. Miller}, \enquote{Computational methods
  of neutron transport,}  (1984).

\bibitem{miften1993quasi}
\textsc{M.~Miften} and \textsc{E.~W. Larsen}, \enquote{The quasi-diffusion
  method for solving transport problems in planar and spherical geometries,}
  \emph{Transport Theory and Statistical Physics}, \textbf{22}, \emph{2-3}, 165
  (1993).

\bibitem{vermaak2021massively}
\textsc{J.~I. Vermaak}, \textsc{J.~C. Ragusa}, \textsc{M.~L. Adams}, and
  \textsc{J.~E. Morel}, \enquote{Massively parallel transport sweeps on meshes
  with cyclic dependencies,} \emph{Journal of Computational Physics},
  \textbf{425}, 109892 (2021).

\bibitem{knoll2011acceleration}
\textsc{D.~Knoll}, \textsc{H.~Park}, and \textsc{C.~Newman},
  \enquote{Acceleration of k-eigenvalue/criticality calculations using the
  Jacobian-free Newton-Krylov method,} \emph{Nuclear Science and Engineering},
  \textbf{167}, \emph{2}, 133 (2011).

\bibitem{park2012nonlinear}
\textsc{H.~Park}, \textsc{D.~Knoll}, and \textsc{C.~Newman}, \enquote{Nonlinear
  acceleration of transport criticality problems,} \emph{Nuclear Science and
  Engineering}, \textbf{172}, \emph{1}, 52 (2012).

\bibitem{petsc-user-ref}
\textsc{S.~Balay}, \textsc{S.~Abhyankar}, \textsc{M.~F. Adams},
  \textsc{S.~Benson}, \textsc{J.~Brown}, \textsc{P.~Brune},
  \textsc{K.~Buschelman}, \textsc{E.~Constantinescu}, \textsc{L.~Dalcin},
  \textsc{A.~Dener}, \textsc{V.~Eijkhout}, \textsc{J.~Faibussowitsch},
  \textsc{W.~D. Gropp}, \textsc{V.~Hapla}, \textsc{T.~Isaac},
  \textsc{P.~Jolivet}, \textsc{D.~Karpeev}, \textsc{D.~Kaushik}, \textsc{M.~G.
  Knepley}, \textsc{F.~Kong}, \textsc{S.~Kruger}, \textsc{D.~A. May},
  \textsc{L.~C. McInnes}, \textsc{R.~T. Mills}, \textsc{L.~Mitchell},
  \textsc{T.~Munson}, \textsc{J.~E. Roman}, \textsc{K.~Rupp},
  \textsc{P.~Sanan}, \textsc{J.~Sarich}, \textsc{B.~F. Smith},
  \textsc{S.~Zampini}, \textsc{H.~Zhang}, \textsc{H.~Zhang}, and
  \textsc{J.~Zhang}, \enquote{{PETSc/TAO} Users Manual,}  ANL-21/39 - Revision
  3.20, Argonne National Laboratory (2023); {10.2172/1968587}.

\bibitem{falgout2002hypre}
\textsc{R.~D. Falgout} and \textsc{U.~M. Yang}, \enquote{hypre: A library of
  high performance preconditioners,} \emph{International Conference on
  computational science}, 632--641, Springer (2002).

\bibitem{geuzaine2009gmsh}
\textsc{C.~Geuzaine} and \textsc{J.-F. Remacle}, \enquote{Gmsh: A 3-D finite
  element mesh generator with built-in pre-and post-processing facilities,}
  \emph{International journal for numerical methods in engineering},
  \textbf{79}, \emph{11}, 1309 (2009).

\bibitem{prince2020diffusion}
\textsc{Z.~Prince}, \textsc{Y.~Wang}, and \textsc{L.~Harbour}, \enquote{A
  diffusion synthetic acceleration approach to k-eigenvalue neutron transport
  using PJFNK,} \emph{Annals of Nuclear Energy}, \textbf{148}, 107714 (2020).

\end{thebibliography}

\end{document}